\documentclass{article}

\usepackage[latin1]{inputenc}
\usepackage{latexsym}
\usepackage{a4wide}
\usepackage{amscd}
\usepackage{graphics}
\usepackage{amsmath}
\usepackage{amssymb}
\usepackage{mathrsfs}
\input xy
\xyoption{all}

\newcommand{\C}{\mathbb C}
\newcommand{\N}{\mathbb N}
\newcommand{\R}{\mathbb R}

\newcommand{\Z}{\mathbb Z}
\newcommand{\sgn}{\mathrm{sign}}

\newcommand{\iEdwards}{\mu_{\scriptstyle{\mathrm{EM}}}}
\newcommand{\ispec}{\mu_{\scriptstyle{\mathrm{Mor}}}}

\newcommand{\Ufrak}{\mathscr U}

\newcommand{\spfl}{\mathrm{sf}}
\newcommand{\coker}{\mathrm{coker}}

\newcommand{\Herm}{\mathrm{Herm}}

\newcommand{\irel}{\mu_{\scriptscriptstyle{\mathrm{rel}}}}

\newcommand{\Id}{\mathrm{Id}}

\newcommand{\Imm}{\mathrm{Im\,}}
\newcommand{\Gl}{\mathrm{Gl}}

\newcommand{\Graph}{\mathrm{Graph\,}}

\newcommand{\diag}{\mathrm{diag\,}}
\numberwithin{equation}{section}
\newcommand{\proof}{{\sl Proof.}\hspace{5pt}}   
\newcommand{\qedhere}{\hfill $\Box$\\}            

\newtheorem{mainthm}{\sc Theorem}           
\newtheorem{thm}{\sc Theorem}[section]      
\newtheorem{lem}[thm]{\sc Lemma}            
\newtheorem{prop}[thm]{\sc  Proposition}     
\newtheorem{defin}[thm]{\sc Definition}      
\newtheorem{rem}[thm]{\sc Remark}       
\newtheorem{notation}[thm]{\sc Notation}    

\title{On a generalized Sturm theorem}
\author{Alessandro Portaluri
\date{May 21, 2009}
\thanks{%
The author was partially supported by MIUR project {\em Variational
Methods and Nonlinear Differential Equations\/}.} }

\begin{document}

\maketitle
\begin{abstract}
Sturm oscillation theorem for second order differential equations
was generalized to systems and higher order equations with positive
leading coefficient by several authors. What we propose here is a
Sturm theorem for indefinite systems with Dirichlet boundary
conditions of the form
\begin{equation*}
 p_{2m}\dfrac{d^{2m}u}{dx^{2m}}+
p_{2m-2}(x)\dfrac{d^{2m-2}u}{dx^{2m-1}}+ \dots
+p_1(x)\dfrac{du}{dx}+ p_0(x) u =0,
\end{equation*}
where  $p_i$ is a smooth path of matrices on the complex
$n$-dimensional vector space $\C^n$ and $p_{2m}$ is the symmetry
represented by $\diag(I_{n-\nu}, -I_\nu)$ for some integer $0\leq\nu
\leq n$.
\end{abstract}


\section{Introduction}\label{sec:intro}

Sturm oscillation theorem deals with differential equations of the
form
\begin{equation}\label{eq:sturmsc}
-(pu')' + qu = 0
\end{equation}
where $'$ denotes  differentiation and where $p$ and $q$ are given
(differentiable) functions with $p>0$. 
Let $V[0,1]$ be the vector space $C^\infty_0([0,1])$ and, for each
$x \in [0,1]$, let us consider the following two quadratic forms on
$V[0,x]$ and $V[0,1]$ respectively given by $q_x(u) := \int_0^x
\Omega[u]dx$ and $q_\lambda(u):=\int_0^1 (\Omega[u]- \lambda
u^2)dx$, where $\Omega[u] = p|u'|^2 + q |u|^2$. Then, the Sturm
oscillation theorem can be restated as follows:
\begin{equation}\label{eq:sturmoscillazionescalare}
\sum_{0 <x<1}\dim\,\ker\,q_{x}(u)\, =\,\sum_{\lambda <0}\dim\,
\ker\,q_{\lambda}(u).
\end{equation}
In this paper we shall study the Dirichlet boundary value problem
for the linear differential equation:
\begin{equation}\label{eq:opdiff}
l(x, D)u:= p_{2m}\dfrac{d^{2m}u}{dx^{2m}}+
p_{2m-2}(x)\dfrac{d^{2m-2}u}{dx^{2m-2}}+ \dots
+p_1(x)\dfrac{du}{dx}+ p_0(x) u =0, \qquad  x \in [0,1]
\end{equation}
where  $p_i$ is a smooth path of matrices on the complex
$n$-dimensional vector space $\C^n$ and $p_{2m}$ is the symmetry
represented by the diagonal block matrix $\diag(I_{n-\nu}, -I_\nu)$
for some integer $0\leq\nu \leq n$, by using the reformulation in
terms of calculus of variations as in \cite{Edw64}. However, we
observe in this respect, that the situation we are dealing with is
completely different and a new approach is needed, since both the
hand-sides of the equation \eqref{eq:sturmoscillazionescalare} are
meaningless. What we propose here is a different definition of both
sides of formula \eqref{eq:sturmoscillazionescalare}. In fact, in
our case, the Morse index is not well-defined and the natural
substitute for the right hand-side will be the {\em spectral flow\/}
of a suitable family of Fredholm Hermitian forms. Furthermore, the
left hand side will be replace by a Maslov-type index obtained by
specifying a suitable
intersection theory in the classical $\Ufrak$-manifolds contest.\\
\noindent It is worth noticing a generalization of the Sturm
oscillation theorem in an indefinite situation recently obtained in
\cite{Zhu05}. Here the author firstly reduces the even order
differential system to a first order Hamiltonian system and then
proves the equality between the spectral flow of a path of unbounded
self-adjoint Fredholm operator and the Maslov index of a suitable
path of Lagrangian subspaces.

\section{Linear preliminaries }\label{sec:basicfacts}

\paragraph{Spectral flow for Fredholm Hermitian forms.}\label{sec:spectralflow}

Let $H$ be a complex separable Hilbert space. A bounded self-adjoint
operator $A:H \to H$ is Fredholm  if $\ker A$ is finite dimensional
its image is closed and $\coker A$ is also finite dimensional. The
topological group $\Gl(H)$ of all automorphisms of $H$ acts
naturally on the space of all self-adjoint Fredholm operators
$\Phi_{S}(H)$  by cogredience sending $A \in \Phi_{S}(H)$ to
$S^*AS$. This induces an action of paths in $\Gl(H)$ on paths in
$\Phi_{S}(H)$. As in the real case, for any path $A \colon [a,b]\to
\Phi_{S}(H)$ there exist a path $M \colon[a,b]\rightarrow \Gl(H)$
and a symmetry $\mathscr J $ (i.e. a bounded linear operator such
that $\mathscr J^2 =\Id$) such that $M^*(t)A(t)M(t)= \mathscr J
+K(t) $ with $K(t)$ compact for each $ t\in [a,b].$ Assuming that
the path $A$ has invertible endpoints and denoting by $\irel$ the
{\em relative Morse index\/}, then the {\em spectral flow of the
path $A$\/} is the integer
\[
\spfl(A,[a,b])\equiv \irel(\mathscr J+K(a),\mathscr J +K(b)),
\]
where $ \mathscr J+K$ is any  compact perturbation of a symmetry
cogredient with $A$. By the properties of the relative Morse index
it is easy to check that this number is well-defined. The spectral
flow $\spfl(A,[a,b])$ is additive and invariant under homotopies
with invertible end points. We refer to \cite{MusPejPor03} for
further details.\\ \noindent A {\em Fredholm Hermitian form\/} on
$H$ is a function $q\colon H\to \R$ such that there exists a bounded
symmetric sesquilinear form $b= b_{q}\colon H\times H\to\C$ with
$q(u)=b(u,u)$  and with $ \ker b $ of finite dimension. We denote by
$\Herm_F(H)$ the set of all Fredholm Hermitian forms. It is possible
to prove that $\Herm_F(H)$ is an open subset (in the operator norm
topology) of $\Herm(H)$ which is stable under perturbations by
weakly continuous Hermitian forms. Moreover a Hermitian form is
called {\em non degenerate\/} if the map $u\rightarrow b_q (u,-) $
is an isomorphism between $H$ and its dual $H^*$. Furthermore a path
of Fredholm Hermitian forms $q\colon [a,b]\rightarrow \Herm_{F}(H)$
with non degenerate end points  $q(a)$ and $q(b)$ will be called
{\em admissible\/}.
\begin{defin}\label{d4}
The spectral flow of an admissible path $q\colon [a,b]\rightarrow
\Herm_{F}(H)$  is given by
\[
\spfl(q,[a,b]):=\spfl(A_{q},[a,b])
\]
where $A_{q(t)}$ is the unique self-adjoint Fredholm operator such
that $\left\langle A_{q(t)}u,u\right\rangle ={q(t)}(u)$, for all
$u\in H$.
\end{defin}
As consequence of the invariance of the spectral flow under
cogredience, it can be proved that this is independent from the
choice of the scalar product. Now given any differentiable path
$q\colon [a,b]\rightarrow \Herm_F(H)$ at the point $t$, then the
derivative $\dot q(t)$ is also a Fredholm Hermitian form. We will
say that a point $t$ is a {\em crossing point\/} if $\ker
b_{q(t)}\neq \{0\} $, and we will say that the crossing point $t$ is
{\em regular\/} if the {\em crossing form\/} $\Gamma(q,t),$ defined
as the restriction of the derivative $\dot q(t)$ to the subspace
$\ker b_{q(t)}$, is nondegenerate. It is possible to prove that
regular crossings are isolated and that the property of having only
regular crossings is generic for paths in $\Herm_F(H)$.
\begin{prop}\label{thm:crossform}
If all crossing points of the path are regular then they are in a
finite number and $\spfl(q,[a,b])= \sum_i  \sgn  \Gamma(q,t_i).$
\end{prop}

\paragraph{The structure of $\Ufrak$-manifolds and the EM-index.}
In this paragraph we will briefly recall some useful facts about the
$\Ufrak$-manifold and we will define the {\em EM-index\/}, by
generalizing the intersection theory proposed by Edwards in
\cite{Edw64}.
\begin{notation}
A (real-valued) Hermitian form on a complex vector space $V$ is a
real valued function $Q$ on $V$ which satisfies:
\begin{enumerate}
\item[(i)] The parallelogram law: $Q[v_1+v_2]-Q[v_1- v_2]=2(Q[v_1] +
Q[v_2])$ for all $v_1, v_2 \in V$;
\item[(ii)] $Q[cv]=|c|^2Q[v]$ for all $c \in \C$ and $v \in V$.
\end{enumerate}
Such a function is of the form $Q[v]=Q[v,v]$ where $Q: V \times V
\to \C$ is a uniquely determined symmetric sesquilinear form.
\end{notation}

\begin{defin}\label{def:complexhypersymplspace}
A {\em superhermitian space\/} is a pair $(S, h)$, where
\begin{itemize}
\item[(i)] $S$ is a complex even dimensional vector space;
\item[(ii)] $h$ is a non degenerate Hermitian form of zero
signature, called {\em superhermitian structure.\/}
\end{itemize}
\end{defin}
We term {\em superlagrangian\/} subspace any subspace $L$ of the
superhermitian space $(S,h)$ of dimension $1/2 \dim S$ on which the
superhermitian structure $h$ vanishes identically and we will refer
with the name of {\em $\Ufrak$-manifold\/} to the set $\Ufrak (S,
h)$ of all superlagrangian subspaces $L$ of $(S,h)$. We observe that
from the topological viewpoint the $\Ufrak$-manifold is homeomorphic
to the unitary group $U(n)$, where $n=\dim S/2.$ We will refer to
\cite[Section 4]{Edw64}, for further details.\\
\noindent Now, given a finite dimensional complex vector space $V$,
let us consider the space $S:=V \oplus V^*$. If $\zeta := (\xi,
\eta)\in S$ and if $\Im$ denotes the imaginary part of a complex
number, $S$ has a naturally associated superhermitian structure
given by $h[\zeta] := \Im \langle \xi, \eta \rangle$ and usually
called {\em standard superhermitian structure.\/}
Every superlagrangian subspace $P_0$ determines a decomposition of
the space of all superlagrangian subspaces as a disjoint union
\[
\Ufrak =\bigcup_{k=0}^n \Ufrak_k(P_0),
\]
where $n =\dim V$ and where, for each $k$, $\Ufrak_k(P_0)$ is the
submanifold of those superlagrangian subspaces which intersect $P_0$
in a subspace of dimension $k$; i.e.
\[
\Ufrak_k(P_0):=\big\{P\in\Ufrak:\dim\,(P\cap P_0)=k\big\}.
\]
We define the following variety
\[
\mathscr S(P_0)= \bigcup_{k=1}^n\Ufrak_k(P_0).
\]
\begin{defin} Given a pair $P_0,P_1 \in \Ufrak(S,h)$
of complementary superlagrangians and identifying $P_0^*$ with $P_1$
via the symmetric sesquilinear form $\langle \cdot, \cdot \rangle$,
we can define the Hermitian form
$\varphi_{P_0,P_1}:\Ufrak_0(P_1)\rightarrow\Herm(P_0)$ as
\[\varphi_{P_0,P_1}(P) \colon P_0 \oplus P_1 \to \R \colon
(u,T_P u) \mapsto   \Im \langle u, -iT_P u \rangle\]
where $T_P:P_0 \to P_1$ is the unique Hermitian operator whose graph
is $P$.
\end{defin}
Otherwise, the Hermitian form $\varphi_{P_0,P_1}$ on $P$ can be
defined in the following way. Let $j : S \to S$ be the unique map
which is the identity on $P_0$ and the multiplication by $-i$ on
$P_1$. Then it can be easily checked that
\begin{equation}\label{eq:cartaedwards}
\varphi_{P_0, P_1}(P)(v)=h[jv] \qquad \forall\, v \in P,
\end{equation}
where $h$ is the standard superhermitian structure.\footnote{%
We observe that, formula \eqref{eq:cartaedwards} is the definition
of the non-trivial invariant $\alpha$ defined by Edwards in
\cite[Section 4]{Edw64} on triples of superlagrangian
subspaces.}\\\noindent Given any $P\in\Ufrak$, it is possible to
define a canonical isomorphism from $T_P\Ufrak$ and $ \Herm(P)$. In
fact, let $P_0,P_1\in\Ufrak$ be a pair of complementary
superlagrangians and let $\varphi_{P_0,P_1}$ be a chart of $\Ufrak$.
Then the differential $d \colon T_P \Ufrak_0(P_1) \to \Herm(P)$ is
the map which send a point $\hat P \in T_P \Ufrak_0(P_1)$ into the
Hermitian form $Q(P, \hat P)$ on $P$ defined as follows. For all
$\varepsilon
>0$ sufficiently small, let us consider the curve $(-\varepsilon,
\varepsilon) \ni t \mapsto P(t)\in \Ufrak$ such that $P(0)=P$ and
$P'(0)= \hat P$. Then
\begin{equation}\label{eq:formuladicrossing}
Q(P, \hat P)(v) :=  \dfrac{d}{dt}\Big|_{t=0} \Im \langle u,
-iT_{P(t)}u\rangle\, = \, \dfrac{d}{dt}\Big|_{t=0} h[jv(t)]
\end{equation}
for $v(t):=(u,T_{P(t)}u)$, $v:=v(0)\in P$ and where $t \mapsto
T_{P(t)}$ is the path of Hermitian operators contained in the domain
of the chart and such that their graphs agrees with the path of
superlagrangian subspaces $t \mapsto P(t)$ in a sufficiently small
neighborhood of $t=0$.  The differential of the chart gives an
isomorphism between $T_P \Ufrak_0(P_1)$ and $\Herm(P)$ and an easy
computation shows that such isomorphism does not depend on the
choice of $P_1$. Summing up, the following result holds.
\begin{prop}\label{thm:atlantecomplesso}
The $\Ufrak$-manifold is a regular algebraic variety of (complex)
dimension $\frac{n(n+1)}{2}$.  Moreover
$\big(\Ufrak_0(P_1),\varphi_{P_0,P_1}\big)$, when $(P_0,P_1)$ runs
in the set of all pairs of complementary superlagrangians form an
atlas of $\Ufrak$. The differential of $\varphi_{P_0, P_1}(P)$ at
$P$ does not depend on the choice of $P_1 \in \Ufrak_0(P_0)$ and
therefore defines a canonical identification of $T_P\Ufrak$ with
$\Herm(P)$.
\end{prop}
Given any differentiable path $p \colon [a,b] \to \Ufrak$, we say
that $p$ has a {\em crossing\/} with $\mathscr S(P_0)$ at the
instant $t=t_0$ if $p(t_0)\in\mathscr S(P_0)$. At each non
transverse crossing time $t_0 \in [a,b]$ we define the {\em crossing
form\/}
\begin{equation}\label{eq:formulagamma}
\Gamma (p, P_0, t_0)\, =Q\big(p(t_0),p'(t_0)\big)\big|_{p(t_0)\cap
P_0}
\end{equation}
and we say that a crossing $t$ is called {\em regular\/} if the
crossing form $\Gamma(p, P_0, t_0)$ is nonsingular.  It is easy to
prove that regular crossings are isolated and therefore on a compact
interval are in a finite number. Moreover if $p(a), p(b) \notin
\mathscr S(P_0)$ then $p$ is said an {\em admissible path\/}.
\begin{rem}\label{rem:suiprodotti}
We observe that given a superhermitian space $(S,h)$ it can be shown
that the pair defined by $(\widetilde S, \widetilde h)$ where $
\widetilde S:=S\oplus S$ and $\widetilde h:= -h\oplus h$ is a
superhermitian space. With this respect and by defining $\widetilde
j$ as $j \oplus j$, then the crossing form can be formally
represented by formulas
\eqref{eq:formuladicrossing}-\eqref{eq:formulagamma} simply by
writing $\widetilde h $ instead of $h$ and $\widetilde j$ insted of
$j$.
\end{rem}
\begin{thm}\label{thm:teoremadiesistenzaeunuicitadiedwardsindex}
Fix $P_0 \in \Ufrak$. Then there exists one and only one map
\[
\iEdwards(\cdot,P_0): C^0([a,b], \Ufrak) \longrightarrow \Z
\]
satisfying the following axioms:
\begin{enumerate}
\item[(i)]({\em Homotopy invariance\/}) If $p_0, p_1
\colon [a,b]\to \Ufrak$ are two homotopic curves of superlagrangian
subspaces with $p(a), p(b) \notin \mathscr S(P_0)$ then they have
the same EM-index.
\item[(ii)]({\em  Catenation\/}) For $a<c<b$, if $p(c) \notin \mathscr
S(P_0)$, then
\[
\iEdwards( p, P_0) \, = \, \iEdwards(p|_{[a,c]},P_0) +
\iEdwards(p|_{[c,b]}, P_0).
\]
\item[(iii)]({\em  Localization\/}) If $P_0:=\C^n
\times\{0\}$ and $p(t):=\Graph\big(H(t)\big)$ where $t\mapsto H(t)$
is an admissible path of Hermitian matrices having only a regular
crossing at $t=t_0$, then we have
\[
\iEdwards(p, P_0) = \frac12 \sgn\, H(t_0+ \varepsilon) - \frac12
\sgn\, H(t_0-\varepsilon),
\]
where $\varepsilon $ is any positive real number.
\end{enumerate}
The integer $\iEdwards(p,P_0)$ is called the {\em Edwards-Maslov
index\/} of $P_0$ or briefly  {\em EM-index\/}.
\end{thm}
\proof Observe that Axioms $(i)-(ii)$ say that the Maslov index is
an homomorphism of the relative homotopy group $\pi_1\big(\Ufrak,
\Ufrak \setminus \mathscr S(P_0)\big)$ into the integers $\Z$. Now,
since $\Ufrak$ is homeomorphic to the unitary group and since
$\Ufrak \backslash \mathscr S(P_0)$ is a cell,  by excision axiom we
have that $\pi_1\big(\Ufrak, \Ufrak \backslash \mathscr
S(P_0)\big)\cong \Z$. The localization axiom will determine this
homomorphism uniquely. It remains only to show that any two curves
of the type described by axiom $(iii)$ are in the same relative
homotopy class. To do so, let $p_1, p_2$ be two such curves. By
using Kato's selection theorem it is not restrictive to assume that
this two curves are of the form
\[p_j(t)=\Delta(-1, -1, \dots, -1, t, 1, 1,\dots, 1) \qquad t \in [-1,1], \
\textrm{and} \ j =1,2,\] where $\Delta$ denotes the diagonal matrix.
Now the thesis follows by the definition of $\Gamma$ at $t_0$ and by
taking into account that
\[
\sgn\,H(t_0 \pm \varepsilon)=\sgn\,H(t_0) \pm \sgn\,\Gamma(p, P_0,
t_0).
\]\qedhere
Since regular crossing are isolated then on a compact interval are
in a finite number and the following result holds.
\begin{prop}\label{thm:regularcrossingsEM}
For an admissible differentiable path $p: [a,b] \to \Ufrak$ having
only regular crossings, we have:
\begin{equation}\label{eq: MaslovdiRobSal}
\iEdwards(p, P_0)= \sum_{t_0\in(a,b)}\sgn\, \Gamma (p, P_0, t_0)
\end{equation}
where we denote by $\sgn$ the signature of a Hermitian form and
where the summation runs over all crossings $t$.
\end{prop}
\proof The proof of this formula follows by local chart computation
obtained by using formulas \eqref{eq:formuladicrossing}-
\eqref{eq:formulagamma} and the localization and concatenation
properties.
\begin{rem}
We observe that in the case of positive definite leading coefficient
this integer coincides with the total intersection index defined by
Edwards in \cite[Section 4]{Edw64}. In fact, it can be proven that
formula \eqref{eq: MaslovdiRobSal} reduces to \cite[Proposition 4.8,
Property (A)]{Edw64}.
\end{rem}


\begin{section}{Variational setting}\label{sec:variationalsetup}

We use the variational approach to \eqref{eq:opdiff} as described in
\cite{Edw64} and we will stick to the notations of that paper. Given
the complex $n$-dimensional Hermitian space $(\C^n, \langle \cdot,
\cdot \rangle)$, for any $m \in \N$ let $\mathscr H^m:=H^m(J,\C^n)$
be the Sobolev space of all $H^m$-maps from the interval $J:=[0,1]$
into $\C^n$.
\begin{defin}
A {\em derivative dependent Hermitian form\/} or a {\em generalized
Sturm form\/}, is the form $ \Omega(x)[u] = \sum_{i,j=0}^m \langle
D^i u(x), \omega_{i,j}(x) D^j u(x)\rangle $, where, each
$\omega_{i,j}$ is a smooth path of $x$-dependent Hermitian matrices
with constant leading coefficient $\omega_{m,m}:=p_{2m}$.
\end{defin}
We observe that a derivative dependent Hermitian form $\Omega(x)[u]$
actually depends on the $m$-jet $u$, $j^mu:=(u(x), \dots,
u^{(m-1)}(x))$ at the point $x$ and it defines a Hermitian form $q
\colon \mathscr H^m \to \R$ by setting $ q(u):= \int_0^1
\Omega(x)[u] dx.$ If $v \in \mathscr H^m$ and $u \in \mathscr
H^{2m}$ then, using integration by parts, the corresponding
sesquilinear form $q(v,u)$ can be written as
\begin{equation}\label{eq:lachiama2.0}
q(v,u)= \int_0^1\langle v(x), l(x,D)u(x) \rangle dx + \phi(v,u)
\end{equation}
where $l(x,D)$ is a  differential operator of the form of
\eqref{eq:opdiff} and $\phi(v,u)$ is a sesquilinear form depending
only on the $(m-1)$-jet, $j^{m} v(x)$ and on the $(2m-1)$-jet
$j^{2m}u(x)$ at the boundary $x=0,1$. Thus,  there exists a unique
linear map $A(x) : \C^{2mn} \to \C^{mn}$ such that
\begin{equation}\label{eq:defindiA}
\phi(v,u)= [\langle j^m v(x), A(x) j^{2m} u(x)\rangle]_{x=0}^1.
\end{equation}
The only specific fact that will be needed is that the entries
$a_{j, 2m-j-1}$ are all equal to $\pm p_{2m}$. \\\noindent Let
$\mathscr H_0^m:=\mathscr H_0^m(J):=\{u \in \mathscr H^m \colon
j^{m} u(0)=0=j^{m} u(1)\}$ and let $q_\Omega$ be the restriction of
the Hermitian form $q$ to $\mathscr H_0^m$. For each $\lambda \in
J$, let us consider the space $\mathscr H_0^m([0, \lambda])$ with
the form $\int_{[0,\lambda]}\Omega(x)dx$. Via the substitution $x
\mapsto \lambda x$, we transfer this form to $\mathscr H_0^m(J)$,
so, we come to the forms $\Omega_\lambda$ and $q_\lambda$  defined
respectively by
\begin{equation}\label{eq:5}
\Omega_\lambda(x)[u] := \sum_{i,j=0}^{m} \langle D^i u(x),
\lambda^{2m-(i+j)}\omega_{i,j}(\lambda x) D^j u(x)\rangle \quad
\textrm{and}\quad q_\lambda(u):= \int_0^1 \Omega_\lambda(x)[u]dx.
\end{equation}
Then $\lambda \mapsto q_\lambda$ is a smooth path of Hermitian forms
acting on $\mathscr H_0^m$ with $q_1=q_{\Omega}$ and with $q_0(u)=
\int_J \langle p_{2m} D^m u, D^m u \rangle dx$. Using integration by
parts, the sesquilinear Hermitian form $q_\lambda(v,u)$ can be
written as
\[
q_\lambda(v,u)=\int_0^1 \langle v(x), l_\lambda(x,D) u(x) \rangle dx
+ \phi_\lambda(v,u),
\]
where
\[
\phi_\lambda(v,u):=[\langle j^m v(x),
A_\lambda(x)j^{2m}u(x)]_{x=0}^1 \ \  \textrm{and}\ \
l_{\lambda}(x,D)
=p_{2m}\frac{d^{2m}}{dx^{2m}}+\sum_{k=0}^{2m-1}p_k(\lambda\,
x)\lambda^{2m-k}\frac{d^k}{dx^k}.
\]
\begin{defin}
A {\em conjugate instant\/} for $q_{\Omega}$ is any point $\lambda
\in(0,1]$ such that $\ker q_\lambda \not=\{0\}$.
\end{defin}

\noindent Let $C_\lambda$ be the path of bounded self-adjoint
Fredholm operators associated to $q_\lambda$  via the Riesz
representation theorem.
\begin{lem}
The Hermitian form $q_0$ is non degenerate. Moreover each
$q_\lambda$ is a Fredholm Hermitian form. (i.e. $C_\lambda$ is a
Fredholm operator). In particular $\dim \ker q_\lambda <+\infty$ and
$q_\lambda$ is non degenerate if and only if $\ker q_\lambda =
\{0\}$.
\end{lem}
\proof That the operator $C_0$ is an isomorphism can be proven
exactly as in \cite[Proposition 3.1]{MusPejPor03}. On the other hand
each $q_\lambda$ is a weakly continuous perturbation of $q_0$ since
it differs from $q_0$ only by derivatives of $u$ of order less than
$m$. Therefore $C_\lambda- C_0$ is compact for all $\lambda \in J$
and hence $C_\lambda$ is Fredholm of index $0$. The last assertion
follows from this.\qedhere

When the form $q_\Omega$ is non degenerate, i.e. when $1$ is not a
conjugate instant, according to the definitions and notation of
section \ref{sec:basicfacts}, we introduce the following definition.
\begin{defin}
The {\em (regularized) Morse index\/} of $q_\Omega$ is defined by
\begin{equation}\label{def:ispettrale}
\ispec(q_\Omega):= \spfl (q_\lambda, J),
\end{equation}
where $\spfl$ denotes the {\em spectral flow of the path
$q_\lambda$\/} i.e., the number of positive eigenvalues of
$C_\lambda$ at $\lambda=0$ which become negative at $\lambda=1$
minus the number of negative eigenvalues of $C_\lambda$ which become
positive. (See, for instance \cite{MusPejPor03}, for a more detailed
exposition).
\end{defin}

\paragraph{Solution space and EM-index of the boundary value problem.}
\begin{defin} Let
$\Omega$ be a derivative dependent Hermitian form on $\mathscr H^m$.
$ u \in \mathscr H^m$ will be called a {\em solution\/} of $\Omega$
if it is orthogonal with respect to $q$ to $\mathscr H_0^m$.
\end{defin}
If $\Sigma$ is the set of all solutions of $\Omega$, by general
facts on ODE, it can be proved that $\Sigma$ is a subspace of
$\mathscr H^m$ of dimension $2mn$ and
\[
h[u]:= \Im \langle j^m u(x), A(x)j^{2m}u(x)\rangle
\]
is {\em independent of the choice of $x$\/}.(See \cite[Section
1]{Edw64}, for further details). Furthermore, $h$ it is immediately
seen to be a non degenerate Hermitian form on $\Sigma$. To prove
this fact, we first introduce the map $A^\#$, as follows:
\[
A^\#(x): \Sigma \longrightarrow \C^{nm} \oplus\C^{nm}, \quad
\textrm{by}\quad A^\#(x)[u]:= (j^m u(x), A(x)j^{2m}u(x) ).
\]
Let $A(x)=[a_{jk}(x)]_{j,k}$. From the equalities $a_{j, 2m-j-1}(x)=
\pm p_{2m}$ and $a_{j,k}(x)=0$ for $|j+k|\geq 2m$ we have that the
matrix $A^\#(x)$ is non singular, hence $A^\#(x)$ is $1-1$ and onto.
Now the non degeneracy of $h$ follows from the fact that the
Hermitian form $(v,w)\mapsto \Im \langle v, w\rangle$ on
$\C^{nm}\oplus \C^{nm}$ is non degenerate.
\begin{defin}
By the solution space of $\Omega$ we mean the pair $(\Sigma, h)$
consisting of the solutions space $\Sigma$ and the non degenerate
Hermitian form $h$ defined on them by
\[
h[u]:= \Im \langle j^m u(x), A(x)j^{2m}u(x)\rangle.
\]
\end{defin}
Before proceeding further we observe that the point $\lambda_0 \in
(0,1)$ is a {\em conjugate point\/} if there exists a non trivial
solution $u$ of the Dirichlet boundary value problem
\begin{equation}\label{eq:bounvalueprobgenturm}
\left\{\begin{array}{ll}
l(x,D)u(x)=  0,  \qquad \forall\   x \in [0,\lambda_0] \\
j^m u(0)=0=j^m u(\lambda_0).
\end{array}\right.
\end{equation}
Now let $\C^{4mn}=(\C^{mn})^4$, and let $\widetilde h$ be the
Hermitian form given by $\widetilde h(v_1, w_1, v_2, w_2):= -\Im
\langle v_1, w_1 \rangle + \Im \langle v_2, w_2\rangle$. By an easy
calculation obtained by taking into account that superhermitian
structure $h[u]$ is independent of $x$, it follows that the image of
the map
\[
\widetilde A^\#:=A^\#(0) \oplus A^\#(1): \Sigma \to \C^{4mn}: u
\mapsto \big(A^\#u(0), A^\#u(1)\big)
\]
is an element of $\Ufrak(\C^{4mn}, \widetilde h)$. It is a
well-known fact that conjugate points cannot accumulate at $0$ and
thus we can find an $\varepsilon >0$ such that there are no
conjugate points in the interval $[0, \varepsilon]$.\\ \noindent Now
denoting by $a: [\varepsilon, 1] \to \Ufrak(\C^{4mn}, \widetilde h)$
the path defined by $a(\lambda):= \Imm \widetilde A_\lambda^\#$, its
EM-index is well-defined and independent on $\varepsilon$. Thus we
are entitled to give the following.
\begin{defin}\label{def:revisitedEdwardssolution1}
Let $P_0:=\{0\}\oplus \C^{mn}\oplus \{0\}\oplus \C^{mn}$. We define
the {\em EM-index\/} of $\Omega$ as the integer given by
\[
\iEdwards(\Omega):= \iEdwards \big(a\vert_{[\varepsilon,1]},
P_0\big).
\]
\end{defin}
\end{section}



\begin{section}{The main result}\label{sec:mainresults}

\begin{mainthm}\label{thm:genmorseindex}{\em(Generalized Sturm
oscillation theorem).\/} Under notations above, we have:
\[
\iEdwards(u)=\ispec(q_\Omega).
\]
\end{mainthm}
\proof We split the proof into some steps. \\ \noindent {\em The
result holds for regular paths.\/} Let $q$ be the path of Fredholm
Hermitian forms defined by
\begin{equation}\label{eq:5bis}
q_\lambda(u):= \int_0^1 \sum_{i,j=0}^{m} \langle D^i u(x),
\lambda^{2m-(i+j)}\omega_{i,j}(\lambda x) D^j u(x)\rangle dx.
\end{equation}
In order to prove the thesis, as consequence of propositions
\ref{thm:crossform} and \ref{thm:regularcrossingsEM}, it is enough
to show that at each crossing point $\lambda_0$, we have
\[
\sgn \,\Gamma(q, \lambda_0)= \sgn\, \Gamma (a, P_0, \lambda_0).
\]
Now let $\lambda_0$ be a crossing point and let us denote by $\cdot
$ the derivative with respect to $\lambda$. Thus we have
\begin{eqnarray*}
\dot q_{\lambda_0}
(v,u)&=&\dfrac{d}{d\lambda}\Big\vert_{\lambda=\lambda_0}\int_0^1
\Omega_\lambda(x) [v,u] dx=
\dfrac{d}{d\lambda}\Big\vert_{\lambda=\lambda_0}\Big[\int_0^1
\langle v(x), l_\lambda(x,D) u(x) \rangle dx + \phi_\lambda(v,u)\Big]\\
&=& \int_0^1 \langle v(x), \dot l_{\lambda_0}(x,D) u(x)\rangle +
\dot\phi_{\lambda_0}(v,u).
\end{eqnarray*}
We set $S(x,D):=\sum_{k=0}^{2m-1}p_k(x)\frac{d^k}{dx^k}$ and
$S_\lambda(x,D)=\lambda^{2m-k} S(\lambda \, x)$. Thus
$l_\lambda(x,D)$ can  be written as $p_{2m}\frac{d^{2m}}{dx^{2m}}+
S_\lambda(x,D)$ and the following result holds.
\begin{lem}
If $u$ is a solution of $l(x,D)u=0$ then
$u_s(x):=u\big(\frac{s}{\lambda}x\big)$ is a solution of
\[
p_{2m}  \frac{d^{2m}}{dx^{2m}}u_s(x) + S_s(x,D)u_s(x)=0.
\]
\end{lem}
\proof It  follows by a straightforward calculations.\qedhere\\
Therefore for any $s \in (0,1]$ the function $u_s(x)$ solves the
Cauchy problem
\begin{equation}\label{cp}
\left\{ \begin{array}{ll} p_{2m}  \frac{d^{2m}}{dx^{2m}}u_s(x) +
S_s(x,D)u_s (x)\equiv
l_s(x,D)u_s(x) = 0\\
u_s (0)=0, u_s' (0) = c_1 u'(0) , \dots , D^{2m-1}u_s (0) =
c_{2m-1}u^{(2m-1)}(0)
\end{array}\right.
\end{equation}
where for each $j=1, \dots, 2m-1$, $c_j= \frac{s^j}{\lambda_0^j}$.
Differentiating the Cauchy problem  \eqref{cp} with respect to $s$
and evaluating at $s=\lambda_0$, we get
\begin{equation}\label{eq:cpinlambda}
\left\{\begin{array}{ll} l_{\lambda_0}(x,D){\dot u}_{\lambda_0} (x)
+ \dot
S_{\lambda_0} (x,D)u_{\lambda_0} (x)=0 \\
\dot u_{\lambda_0} (0)=\dots =\dot u^{(2m-1)}_{\lambda_0} (0) = 0.
\end{array}\right.
\end{equation}
If  $u \in \ker q_{\lambda_0}$, performing integration by parts and
observing that $u_{\lambda_0}(\cdot)=u(\cdot)$, as consequence of
equation \eqref{eq:cpinlambda}, we have
\begin{eqnarray*}
\dot q_{\lambda_0}(u)&=&\int_0^1 \langle u_{\lambda}(x), \dot
l_{\lambda_0}(x,D) u_{\lambda}(x)\rangle +
\dot\phi_{\lambda}(u_{\lambda})= -\int_0^1 \langle u_{\lambda}(x),
l_{\lambda_0}(x,D)\dot u_{\lambda} (x) \rangle dx+
\dot\phi_{\lambda_0}(u_{\lambda})
=\\
&=&-\int_0^1\langle l_{\lambda_0}(x,D) u_{\lambda}(x), \dot
u_{\lambda}(x)\rangle dx + \dot\phi_{\lambda_0}(u_{\lambda})=
\dot\phi_{\lambda_0}(u_{\lambda})=\dot\phi_{\lambda_0}(u).
\end{eqnarray*}
Moreover
\begin{eqnarray*}
\dot\phi_{\lambda_0}(u)&=&\dfrac{d}{d\lambda}\big\vert_{\lambda=\lambda_0}
\Big\{\big[\langle j^m u(x), A_\lambda(x)
j^{2m}u(x)\rangle\big]_{x=0}^1\Big\}= \big[\langle j^m  u(x), \dot
A_{\lambda_0}(x) j^{2m}u(x)\rangle\big]_{x=0}^1=\\ &=&   - \langle
j^m u(0), \dot A_{\lambda_0}(0) j^{2m}u(0)\rangle +\langle j^m u(1),
\dot A_{\lambda_0}(1) j^{2m}u(1)\rangle.
\end{eqnarray*}
Since $\dot q_{\lambda_0}$ is a Hermitian form, in particular it is
a real-valued function; thus we have
\begin{eqnarray*}
\dot q_{\lambda_0}(u)&=& \Re\dot q_{\lambda_0}(u)= \Re
\dot\phi_{\lambda_0}(u)=\Re\big[\langle j^m u(x), \dot
A_{\lambda_0}(x) j^{2m}u(x)\rangle\big]_{x=0}^1=\\ &=& -\Re \langle
j^m u(0), \dot A_{\lambda_0}(0) j^{2m}u(0)\rangle +\Re\langle j^m
u(1), \dot A_{\lambda_0}(1) j^{2m}u(1)\rangle.
\end{eqnarray*}
Since $\Re \langle u,v \rangle = \Im \langle u, -i v\rangle$, we can
conclude that
\begin{eqnarray*}
\dot q_{\lambda_0}(u)&=& -\Re \langle j^m u(0), \dot
A_{\lambda_0}(0) j^{2m}u(0)\rangle +\Re\langle j^m u(1), \dot
A_{\lambda_0}(1) j^{2m}u(1)\rangle=\\
&=&-\Im \langle j^m u(0), -i\dot A_{\lambda_0}(0) j^{2m}u(0)\rangle
+\Im\langle j^m u(1), -i\dot A_{\lambda_0}(1) j^{2m}u(1)\rangle=\\
&=& -h\big[j\dot A_{\lambda_0}^\#(0)[u]\big]+ h\big[j\dot
A_{\lambda_0}^\#(1)[u]\big]= \widetilde h[\widetilde j \dot u^\#
(\lambda_0)]=\\
&=& \Gamma (a, P_0, \lambda_0)(u)
\end{eqnarray*}
where, for $k=0,1$, we denoted by $\dot A_{\lambda_0}^\#(k)[u]$ the
pair $\big(j^m u(k), \dot A_{\lambda_0}(k) j^{2m}u(k)\big)$, by $
\dot u^\#(\lambda_0)$ the element $\big(j^m u(0), \dot
A_{\lambda_0}(0) j^{2m}u(0), j^m u(1), \dot A_{\lambda_0}(1)
j^{2m}u(1)\big)\in \dot a(\lambda_0)$ and where the last equality
follows by remark \ref{rem:suiprodotti}. The above calculations
shown that regular crossings of $q$ correspond to regular crossings
of $a$. Furthermore, the crossing forms at each crossing point
associated to the path of Fredholm Hermitian forms and to the path
of superlagrangian subspaces are the same and therefore their
signatures coincide; in symbols we have
\[
\sgn \,\Gamma(q, \lambda_0)= \sgn\, \Gamma (a, P_0, \lambda_0).
\]
Now the conclusion of the first step follows by
the previous calculations and by summing over all crossings.\\
\noindent {\em Second step. The general case.\/} In order to
conclude remains to show that it is possible to extend the above
calculation to general paths having not only regular crossings. For
each $\lambda \in [0,1]$ let us consider the closed unbounded
Fredholm operator $A_\lambda$ on $L^2(J, \C^{mn})$ with domain
$\mathscr D( A_\lambda) =\{u \in \mathscr H^{2m}:j^{m-1}
u(0)=0=j^{m-1} u(1) \}$ defined by $ A_\lambda u :=
l_\lambda(x,D)u$. By applying a perturbation argument proven in
\cite[Theorem 4.22]{RobSal95} to the path of operators $A_\lambda$,
we can find a $\delta>0$ such that $A_\lambda^\delta:= A_\lambda +
\delta \Id$ is a path of self-adjoint Fredholm operators with only
regular crossing points. Let $q_\lambda^\delta(u)$ be the Hermitian
form on $\mathscr H_0^m$ given by  $ q_\lambda^\delta(u):=\langle u,
A_\lambda^\delta (u)\rangle_{L^2}+ \frac12 \delta\|u\|_{L^2}^2$. By
this choice of $\delta $ and by applying the first step to the
perturbed path $q_\lambda^\delta$, we conclude the proof.\qedhere
\end{section}


\end{document}